\documentclass[12pt]{amsart}
\usepackage{graphicx}
\usepackage{mathrsfs}

\newtheorem{thm}{Theorem}[section]

\newtheorem{lem}[thm]{Lemma}

\theoremstyle{definition}

\theoremstyle{remark}

\numberwithin{equation}{section}

\DeclareMathOperator*{\cov}{cov} \DeclareMathOperator*{\var}{var}
\usepackage{latexsym,amsfonts,amsmath}
\begin{document}

\title[Branching populations before extinction]{Asymptotic behavior of
branching populations before extinction}%
\author{Vyacheslav M. Abramov}%
\address{School of Mathematical Sciences, Building 28M, Monash University, Clayton
Campus, Clayton, VIC 3800, Australia}%
\email{Vyacheslav.Abramov@sci.monash.edu.au}%


\begin{abstract}
Under the assumption that the initial population size of a
Galton-Watson branching process increases to infinity, the paper
studies asymptotic behavior of the population size  before
extinction. More specifically, we establish asymptotic properties
of the conditional moments (which are exactly defined in the
paper).
\end{abstract}
\maketitle
\section{\bf Introduction and the main result}\label{sec-1}
\noindent We consider a Galton-Watson branching process
$\{X_n\}_{n\geq0}$,
\begin{equation}\label{eq0}
X_{n+1}=\sum_{j=1}^{X_n}\xi_{n,j},
\end{equation}
where $X_n$ denotes the number of offspring in the $n$th
generation for a population starting from $K$ offsprings, $X_0=K$,
and throughout the paper the initial size of population $K$ is
assumed to be a large value. Such a type of branching process can
be a model of real population of animals, insects etc., and the
main results of our study can have applications to analysis of
real populations arising in biology (e.g. Jagers \cite{Ja 1975},
Haccou, Jagers and Vatutin \cite{HJV 2005}, Jagers and Klebaner
\cite{Ja-Kl 2003}). For other study of branching processes with a
large initial population size see also Borovkov \cite{Bor 1991},
Klebaner \cite{Kl 1993}, Klebaner and Liptser \cite{Kl-Lip in
preparation}.

The study of branching populations before extinction has been
initiated by Jagers \cite{Ja before ex} and then resulted in
papers of Jagers, Klebaner and Sagitov \cite{Ja 2005} and \cite{Ja
2007}. The approach of these papers \cite{Ja 2005} and \cite{Ja
2007} is based on analytic techniques for studying paths to
extinction with the following analysis of asymptotic behavior of
these paths.

The present paper presents an alternative way to study asymptotic
behavior of large populations before extinction, and the approach
of the present paper is based on diffusion approximations of the
original branching process with large initial population as well
as a series of auxiliary processes. Those diffusion approximations
are then used to study asymptotic behavior of conditional moments
of a population size before extinction as it is explained below.
The approach of our paper (including diffusion approximations,
asymptotic expansions and sample path techniques) remains correct
for much wider classes of branching processes than that
traditional branching process and includes for instance bisexual
Galton-Watson branching processes \cite{Bruss 1984}, \cite{Daley
1968} and different type of controlled $\phi$-branching
Galton-Watson processes (which need not be subcritical, as it is
assumed in the paper later). The $\phi$-branching processes have
been introduced by Sevastyanov and Zubkov \cite{Sev-Zub 1974}, and
intensively studied in many papers (e.g. Bruss \cite{Bruss 1978},
\cite{Bruss 1980}, Gonz\'alez, Molina and Del Puerto \cite{GMDP1},
\cite{GMDP2}, \cite{GMDP3}, Yanev \cite{Yan 1975}, Zubkov
\cite{Zub 1974}).

The main results of the present paper are presented in Theorem
\ref{th-1} below.

\smallskip

Assume that $\xi_{n,j}$, ($n$=1,2,\ldots; $j$=1,2,\ldots) have the
same distribution for all $n$ and $j$ and are mutually
independent, and there exists the second moment
$\mathrm{E}\xi_{n,j}^2<\infty$. Denoting
$\frak{m}=\mathrm{E}\xi_{n,j}$ and $\frak{S}^2=\var(\xi_{n,j})$,
assume that $\frak{m}< 1$. Under this last assumption the
extinction time of the branching process always exists with
probability 1. Let $\tau=\tau_K$ be that moment of extinction. The
random variable $\tau_K$ is a stopping time associated with the
sequence $\{X_n\}_{n\geq0}$. We assume that the family of all
stopping times $\{\tau_K\}$ (for different values $K$) is defined
on a filtered probability space $\{\Omega, \mathscr{F}_0,
\mathbf{F}=(\mathscr{F}_{0,K}), \mathrm{P}\}$,
$\mathscr{F}_{0,K}\subset\mathscr{F}_{0,K+1}\subset\ldots\subset\mathscr{F}_0$.
(The meaning of the index 0 will be clear later.)

The paper studies asymptotic behavior of the branching population
before extinction as $K$ increases to infinity, and the main
result of our study, formulated below, as well as the analysis of
the paper, use the notation $\asymp$ for asymptotic equivalence
between two main parts of expansion. The notation is used in order
to reduce irrelevant background explanations and to avoid multiple
using of $\lim$ in different senses or expansions with remainder,
where it is not significant. For example, relations
\eqref{eq-th-1} and \eqref{eq-th-3} (see formulation of the
theorem below) should be read as follows: \textit{For any
sufficiently small positive $\epsilon$ and $\delta$ there exists a
large integer $K$ such that}
\begin{equation*}
\begin{aligned}
&\mathrm{P}\left\{(1-\delta)X_{\lfloor
u_2\tau_K\rfloor}^l\mathrm{E}\frak{m}^ {l(\lfloor
u_1\tau_{K}\rfloor-\lfloor u_2\tau_{K}\rfloor)}\right.\\
&\ \ \ \left.\leq\mathrm{E}(X_{\lfloor
u_1\tau_K\rfloor}^l~|~X_{\lfloor u_2\tau_K\rfloor}) \leq(1+\delta)
X_{\lfloor u_2\tau_K\rfloor}^l\mathrm{E}\frak{m}^ {l(\lfloor
u_1\tau_{K}\rfloor-\lfloor u_2\tau_{K}\rfloor)}\right\}>1-\epsilon,\\
&\mathrm{P}\left\{(1-\delta)X_{\lfloor
u_1\tau_K\rfloor}^l\mathrm{E}\frak{m}^ {l(\lfloor
u_2\tau_{K}\rfloor-\lfloor u_1\tau_{K}\rfloor)}\right.\\
&\ \ \ \left.\leq\mathrm{E}(X_{\lfloor
u_2\tau_K\rfloor}^l~|~X_{\lfloor u_1\tau_K\rfloor}) \leq(1+\delta)
X_{\lfloor u_1\tau_K\rfloor}^l\mathrm{E}\frak{m}^
{l(\lfloor u_2\tau_{K}\rfloor-\lfloor u_1\tau_{K}\rfloor)}\right\}>1-\epsilon,\\
\end{aligned}
\end{equation*}
and
\begin{equation*}
\begin{aligned}
&\mathrm{P}\left\{(1-\delta)K^l\frak{m}^{l\lfloor
u_1\tau_K\rfloor}\right.\\ &\ \ \ \left.\leq\mathrm{E}(X_{\lfloor
u_1\tau_K\rfloor}^l~|~\tau_K)\leq(1+\delta)K^l\frak{m}^{l\lfloor
u_1\tau_K\rfloor}\right\}>1-\epsilon.
\end{aligned}
\end{equation*}

In the places where it is required and looks more profitable (e.g.
Section \ref{sec-5}), the explicit form of asymptotic expansion
with remainder is used nevertheless.

\begin{thm}\label{th-1}
Let $0<u_1<u_2<1$ be two real numbers. Then, as $K\to\infty$,
\begin{equation}\label{eq-th-1}
\begin{aligned}
&\mathrm{E}\{X_{\lfloor u_1\tau_K\rfloor}^l~|~X_{\lfloor
u_2\tau_K\rfloor}\}\asymp X_{\lfloor
u_2\tau_K\rfloor}^l\mathrm{E}\frak{m}^
{l(\lfloor u_1\tau_{K}\rfloor-\lfloor u_2\tau_{K}\rfloor)},\\
&\mathrm{E}\{X_{\lfloor u_2\tau_K\rfloor}^l~|~X_{\lfloor
u_1\tau_K\rfloor}\}\asymp X_{\lfloor
u_1\tau_K\rfloor}^l\mathrm{E}\frak{m}^
{l(\lfloor u_2\tau_{K}\rfloor-\lfloor u_1\tau_{K}\rfloor)},\\
\end{aligned}
\end{equation}
and
\begin{equation}\label{eq-th-3}
\mathrm{E}\{X_{\lfloor u_1\tau_K\rfloor}^l~|~\tau_K\}\asymp
K^l\frak{m}^{l\lfloor u_1\tau_K\rfloor},
\end{equation}
where $\lfloor z\rfloor$ is the notation for the integer part of
$z$.
 As
$K\to\infty$, $\frac{\tau_K}{\log K}$ converges in probability to
the constant $c=-\frac{1}{\log\frak{m}}$.
\end{thm}

The proof of the main result is based on the following lemma.

\begin{lem}\label{lem1}
For any finite-dimensional vector
$\{X_{i_1},X_{i_2},\ldots,X_{i_n}\}$, $1\leq
i_1<i_2<\ldots<i_n<\infty$,
\begin{equation*}
\begin{aligned}
&\lim_{K\to\infty}\mathrm{P}\left\{\frac{X_{i_1}-K\frak{m}^{i_1}}{\frak{S}\sqrt{K}}\leq
x_1, \frac{X_{i_2}-K\frak{m}^{i_2}}{\frak{S}\sqrt{K}} \leq
x_2,\ldots,\frac{X_{i_n}-K\frak{m}^{i_n}}{\frak{S}\sqrt{K}}\leq
x_n \right\}\\
&=\mathrm{P}\{\theta_{i_1}\leq x_1, \theta_{i_2}\leq x_2,\ldots,
\theta_{i_n}\leq x_n\},
\end{aligned}
\end{equation*}
where $\{\theta_1,\theta_2,\ldots\}$ is a Gaussian sequence with
$\mathrm{E}\theta_j=0$ and $\cov(\theta_{j}$, $\theta_{j+n}) =
\frak{m}^n\var(\theta_j) + n\frak{m}^{j+n-1}$, $\var(\theta_{j+1})
= \frak{m}^2\var(\theta_j) + \frak{m}^j$, \ $\var(\theta_1)=1$.
\end{lem}

Lemma \ref{lem1} is known from the literature, and its proof can
be found in Klebaner and Nerman \cite{Kl-Ner}. For the purpose of
the present paper we, however, need in an alternative proof of
this lemma, which follows from the asymptotic expansions presented
here. Furthermore, the proof of Theorem \ref{th-1} requires the
intermediate asymptotic expansions obtained in the proof of Lemma
\ref{lem1} rather than the statement of Lemma 1.2 itself.

In this paper, simple asymptotic representations for all
conditional moments before extinction are obtained. The most
significant consequence of this analysis is a so-called
\textit{invariance} property of the conditional expectations. This
property is discussed in Section \ref{sec-5}.

The main idea of the method is as follows. The random sequence
$\{X_n\}_{n\geq0}$ is approximated by appropriate random sequences
$\{Y_n^{(a)}\}_{n\geq0}$ ($a\leq1$), as $a$ tends to zero. For
each fixed $a$ we define stopping times $\tau_{a,K}$ (for
different values $K$) associated with the process $Y_n^{(a)}$.
$\tau_{a,K}$ is assumed to be measurable with respect to the
$\sigma$-field $\mathscr{F}_{a,K}\subset\mathscr{F}_{a}$, where
$\mathscr{F}_{a}=\cup_{K\geq 1}\mathscr{F}_{a,K}$, and
$\mathscr{F}_a\subset\mathscr{F}_0$. For that fixed $a$ the
sequence $\tau_{a,K}$ converges (in definite sense) to $\ell({a})$
as $K \to\infty$ (the details are given in the paper). Then
knowledge of the behavior of $Y_{u\tau_{a,K}}^{(a)}$, $0<u<1$, for
which we have the corresponding relationship, enables us to study
the behavior of its limit as $a$ tends to zero. This limit is just
$X_{u\tau_{K}}$, $0<u<1$. Other assumptions associated with
definition of $X_n$ and that of the associated processes
$X_n^{(a)}$, $Y_n^{(a)}$ and other processes are given in the next
section.

The rest of the paper is organized as follows. In Section
\ref{sec-2} we introduce the auxiliary stochastic sequences
$X_n^{(a)}$ and $Y_n^{(a)}$ and the stopping times associated with
these sequences. The elementary properties of these random objects
are studied. In Section \ref{sec-3} we continue to study the
properties of the sequences $X_n^{(a)}$ and $Y_n^{(a)}$.
Specifically, it is shown that these sequences are upper and lower
bounds for the branching process $X_n$, and these bounds are tight
as $a\to 0$. These properties are then used in order to prove the
convergence results in the next sections. In Section \ref{sec-4}
we derive asymptotic expansions and prove the convergence lemma to
the Gaussian process, the parameters of which are explicitly
defined in the formulation of Lemma \ref{lem1}. In Section
\ref{sec-proof} we prove Theorem \ref{th-1}. Last Section
\ref{sec-5} discusses application of the main results of this
study and establishes the invariance property.

\section{\bf Stopping times and auxiliary processes associated with the Galton-Watson
process}\label{sec-2} \noindent In this section we approach the
stopping time $\tau_{K}$, the extinction moment, by introducing a
parametric family of stopping times $\{\tau_{a,K}\}$, depending on
the two parameters $a$ and $K$. Specifically, for any real $a$,
$0\leq a<1$ and integer $K$
\begin{equation}\label{eq1}
\tau_{a,K}=\inf\{l: X_{l}\leq \lfloor aK\rfloor\},
\end{equation}
where $\lfloor aK\rfloor$ is the integer part of $aK$. The
stopping time $\tau_{a,K}$ as well as the associated with these
parameters $a$ and $K$ other corresponding random variables
defined below are assumed to be measurable with respect
$\mathscr{F}_{a,K}\subset \mathcal{F}_K$, and for two different
values $a_1$ and $a_2$, $0\leq a_2<a_1<1$, we have
$\mathscr{F}_{a_1,K}\subset\mathscr{F}_{a_2,K}$. If $a<1$ is fixed
and $K_1$, $K_2$ are distinct, $K_1<K_2<\infty$, then we have
$\mathscr{F}_{a,K_1}\subset\mathscr{F}_{a,K_2}$. Then the
two-parametric family of $\sigma$-fields $\{\mathscr{F}_{a,K}\}$
is increasing in the following sense. For any $0\leq a_2\leq
a_1<1$ and integer $K_1\leq K_2<\infty$ we have
$\mathscr{F}_{a_1,K_1}\subseteq\mathscr{F}_{a_2,K_2}$.

In accordance with this family of stopping times \eqref{eq1},
consider a family of processes $X_{j,K}^{(a)}=X_j^{(a)}$
satisfying the recurrence relation (for notational convenience the
additional index $K$ is not provided):
\begin{equation}\label{eq2}
X_{n+1}^{(a)}=\max\left\{\lfloor aK\rfloor,
\sum_{j=1}^{X_n^{(a)}}\xi_{n,j}\right\}, \ \ X_0^{(a)}=K.
\end{equation}

The processes $X_{j,K}^{(a)}$ are assumed to be adapted with
respect to the $\sigma$-fields $\mathscr{F}_{a,K}$. In addition,
the processes $X_{j,K}^{(a)}$ are assumed to be measurable with
respect to the the wider $\sigma$-field $\mathscr{F}_{0}$.
Specifically, if there are two processes $X_{j,K}^{(a_1)}$ and
$X_{j,K}^{(a_2)}$ with different $a_1$ and $a_2$, say $0\leq
a_2\leq a_1<1$, then both of these processes $X_{j,K}^{(a_1)}$ and
$X_{j,K}^{(a_2)}$ are measurable with respect to the
$\sigma$-field $\mathscr{F}_{0,K}$, and, of course, with respect
to the $\sigma$-field $\mathscr{F}_{a_2,K}$. All of these
processes with different $a$ are defined due to representation
\eqref{eq2}. This means that the processes $X_{j,K}^{(a)}$ are
actually defined after their stopping times as well. For different
$a_1$ and $a_2$ ($0\leq a_2\leq a_1<1$) the processes
$X_{j,K}^{(a_1)}$ and $X_{j,K}^{(a_2)}$ are `coupled' until the
stopping time $\tau_{a_1,K}$, i.e. until that time instant their
sample paths coincide, but after the time instant $\tau_{a_1,K}$
these processes are decoupled i.e. their paths become different.
But the coupling arguments can be used nevertheless: after the
time instant $\tau_{a_1,K}$ with the aid of Kalmykov's theorem
\cite{Kalmykov} we have $X_{j,K}^{(a_1)}\geq_{st}X_{j,K}^{(a_2)}$,
$j\geq\tau_{a_1,K}$ (see the next section for details).

 Some mathematical details about these processes can be found in
the next section. The similar coupling arguments hold for the
processes $Y_{j,K}^{(a)}$ defined later, which are derivative from
the processes $X_{j,K}^{(a)}$ (the further details can be found in
the next section).

Let us transform \eqref{eq2} by adding and subtracting the term
$\lfloor aK\rfloor$. To this end we use the following elementary
property of numbers: $\max\{a,b\}-a=\max\{0, b-a\}$. Also there is
used the fact that $X_n^{(a)}\geq \lfloor aK\rfloor$ for any $n$.
Then, we have
\begin{equation}\label{eq3}
\begin{aligned}
&\Big(X_{n+1}^{(a)}-\lfloor aK\rfloor\Big)+\lfloor aK\rfloor\\
&=\left(\max\left\{\lfloor aK\rfloor,
\sum_{j=1}^{X_n^{(a)}}\xi_{n,j}\right\}-\lfloor aK\rfloor\right)
+\lfloor aK\rfloor\\
&=\max\left\{0, \sum_{j=1}^{
X_n^{(a)}}\xi_{n,j}-\sum_{j=1}^{\lfloor aK\rfloor}\Big[\xi_{n,j}+
\Big(1-\xi_{n,j}\Big)\Big]\right\}+\lfloor aK\rfloor\\
&=\max\left\{0, \sum_{j=\lfloor aK\rfloor+1}^{
X_n^{(a)}}\xi_{n,j}-\sum_{j=1}^{\lfloor
aK\rfloor}\Big(1-\xi_{n,j}\Big)\right\}+\lfloor aK\rfloor.
\end{aligned}
\end{equation}
Hence, denoting $Y_n^{(a)}=X_n^{(a)}-\lfloor aK\rfloor$ from
\eqref{eq3} we obtain
\begin{equation}\label{eq4}
\begin{aligned}
Y_{n+1}^{(a)}&=\max\left\{0, \sum_{j=\lfloor aK\rfloor+1}^{
X_n^{(a)}}\xi_{n,j}-\sum_{j=1}^{\lfloor
aK\rfloor}\Big(1-\xi_{n,j}\Big)\right\}\\
&=\max\left\{0,
\sum_{j=1}^{Y_n^{(a)}}\xi_{n,j}^\prime-\sum_{j=1}^{\lfloor
aK\rfloor}\Big(1-\xi_{n,j}\Big)
\right\}\\
&=\left(\sum_{j=1}^{Y_n^{(a)}}\xi_{n,j}^\prime
-\sum_{j=1}^{\lfloor aK\rfloor}\Big(1-\xi_{n,j}\Big)\right)\\& \ \
\ \times\mathrm{I}
\left(\sum_{j=1}^{Y_n^{(a)}}\xi_{n,j}^\prime>\sum_{j=1}^{\lfloor
aK\rfloor}\Big(
1-\xi_{n,j}\Big)\right)\\
&=\left(\sum_{j=1}^{Y_n^{(a)}}\xi_{n,j}^\prime-\sum_{j=1}^{\lfloor
aK\rfloor} \Big(1-\xi_{n,j}\Big)\right)I_n^{(a)},
\end{aligned}
\end{equation}
where $\xi_{n,j}^\prime=\xi_{n,j+\lfloor aK\rfloor}$
($\xi_{n,1}^\prime$, $\xi_{n,2}^\prime$,\ldots are independent and
identically distributed random variables having the same
distribution as $\xi_{n,j}$), and
$I_n^{(a)}=I_n^{(a)}(K)=\mathrm{I}
\big\{\sum_{j=1}^{Y_n^{(a)}}\xi_{n,j}^\prime>\sum_{j=1}^{\lfloor
aK\rfloor}(1-\xi_{n,j})\big\}$ is the notation used in
\eqref{eq4}.

Thus we have the new family of processes $Y_n^{(a)}$, which is
assumed, as mentioned before, to be measurable with respect to
$\mathscr{F}_0$ and given on the same probability space
$\{\Omega$, $\mathscr{F}_0$, $\mathrm{P}\}$. Recall that a
stopping time $\tau_{a,K}$ and the sequence
\begin{equation}\label{eq4-1}
\left\{Y_{0}^{(a)}, Y_{1}^{(a)},\ldots\right\}
\end{equation}
are assumed to be adapted with respect to the $\sigma$-field
$\mathscr{F}_{a,K}$, and the family of these $\sigma$-fields
$\{\mathscr{F}_{a,K}\}$ is increasing when $a$ decreases and $K$
increases.

It is known that as $K\to\infty$, $\frac{X_n}{K}$ converges to
$\frak{m}^n$ in probability (see Klebaner and Nerman
\cite{Kl-Ner}). Using this result it is not difficult to prove
that, as $K\to\infty$,

\begin{equation*}
 \frac{Y_n^{(a)}}{K} \ \ \mbox{converges to} \ \max\{0,\frak{m}^n-a\} \ \mbox{in
probability},\leqno(i)
\end{equation*}
\begin{equation*}
\frac{\mathrm{E}Y_n^{(a)}}{K} \ \ \mbox{converges to}\
\max\{0,\frak{m}^n-a\},\leqno(ii)
\end{equation*}
\begin{equation*}
\begin{aligned}
I_n^{(a)}(K) \ \ \mbox{converges to} \ \
\chi_{n+1}&=\chi_{n+1}(\frak{m},a)\\&=
\begin{cases} 1, & \mbox{if}
\ \frak{m}^{n+1}>a,\\
0, & \mbox{otherwise}
\end{cases}
\end{aligned}\leqno(iii)
\end{equation*}

\hskip .85in in probability,

\noindent as well as,
\begin{equation*}
\hskip .40in \tau_{a,K} \ \text{converges in probability to} \
\ell(a)=\min\{l: \frak{m}^l\leq a\}.\leqno(iv)
\end{equation*}
The proof of $(i)$ is postponed to the end of Section \ref{sec-3}.
The proofs of $(ii)-(iv)$ are similar to the proof of $(i)$.

\section{\bf Properties of the sequences
$X_n^{(a)}$ and $Y_n^{(a)}$}\label{sec-3}\noindent The study of
this section  we start from the properties of the random vectors
\eqref{eq4-1}. Let $a_1$, $a_2$ be two numbers, and $0\leq a_2\leq
a_1<1$ and $K$ is fixed. Then, in the suitable probability space
for all events $\omega\in\Omega$ and $n\geq 0$
\begin{equation}\label{eq6}
Y_{n}^{(a_1)}(\omega)\leq Y_n^{(a_2)}(\omega).
\end{equation}
Indeed, consider two random vectors
\begin{equation}\label{eq4-2}
\left\{Y_{0}^{(a_1)}, Y_{1}^{(a_1)},\ldots\right\}
\end{equation}
and
\begin{equation}\label{eq4-3}
\left\{Y_{0}^{(a_2)}, Y_{1}^{(a_2)},\ldots\right\}
\end{equation}
Consider the stopping times $\{\tau_{a_1,K},\mathscr{F}_{a_1,K}\}$
and $\{\tau_{a_2,K},\mathscr{F}_{a_2,K}\}$ associated with the
sequences \eqref{eq4-2} and \eqref{eq4-3}. Since for fixed $K$,
$\mathscr{F}_{a_1,K}\subseteq\mathscr{F}_{a_2,K}$, then
$\tau_{a_1,K}(\omega)\leq\tau_{a_2,K}(\omega)$.

According to the definition of the sequence $X_n^{(a)}$ (see
\eqref{eq2}), on the suitable probability space containing
$\mathscr{F}_{a_2,K}$ we have the correspondence
\begin{equation}\label{eq-3-0}
X_i^{(a_1)}(\omega)=X_i^{(a_2)}(\omega), \ \
i=1,2,\ldots,\tau_{a_1,K}-1,
\end{equation}
and at this stopping time $\tau_{a_1,K}$ we have
$X_{\tau_{a_1,K}}^{(a_1)}(\omega)\geq
X_{\tau_{a_1,K}}^{(a_2)}(\omega)$, and therefore according to
Kalmykov's theorem \cite{Kalmykov}:
\begin{equation*}
X_i^{(a_1)}\geq_{st} X_i^{(a_2)}, \ \
i=\tau_{a_1,K},\tau_{a_1,K}+1,\ldots,\tau_{a_2,K},\ldots,\tau_{0,K}.
\end{equation*}
Therefore, in a suitable probability space
\begin{equation}\label{eq-3-0'}
X_i^{(a_1)}(\omega)\geq X_i^{(a_2)}(\omega), \ \
i=\tau_{a_1,K},\tau_{a_1,K}+1,\ldots,\tau_{a_2,K},\ldots,\tau_{0,K}.
\end{equation}
Thus, we showed
\begin{equation*}
X_i^{(a_1)}(\omega)\geq X_i^{(a_2)}(\omega), \ \
i=1,2,\ldots,\tau_{a_1,K},\tau_{a_1,K}+1,\ldots,\tau_{a_2,K},\ldots,\tau_{0,K}.
\end{equation*}

From this correspondence \eqref{eq-3-0} and \eqref{eq-3-0'}
according to the definition of the sequence $Y_n^{(a)}$ (see
\eqref{eq4}) on the same probability space we have
\begin{equation}\label{eq-3-1}
Y_i^{(a_1)}(\omega)+\lfloor
a_1K\rfloor=Y_i^{(a_2)}(\omega)+\lfloor a_2K\rfloor, \ \
i=1,\ldots,\tau_{a_1,K}-1,
\end{equation}
and therefore up to time $\tau_{a_1,K}-1$ the inequality
$Y_i^{(a_1)}(\omega)\leq Y_i^{(a_2)}(\omega)$ is obvious. At time
instant $\tau_{a_1,K}$ we have
$Y_{\tau_{a_1,K}}^{(a_1)}(\omega)=0$, while
$Y_{\tau_{a_1,K}}^{(a_2)}(\omega)$ is nonnegative in general. The
further behavior of the processes $Y_i^{(a_1)}(\omega)$ and
$Y_i^{(a_2)}(\omega)$ after time $\tau_{a_1,K}$ is specified by
coupling arguments, where the initial inequality
$Y_i^{(a_1)}(\omega)\leq Y_i^{(a_2)}(\omega)$ before the stopping
time $\tau_{a_1,K}$ remains true after this stopping time as well.
If for some $i=i_0$,
$Y_{i_0}^{(a_1)}(\omega)=Y_{i_0}^{(a_2)}(\omega)$ ($=0$), then the
both processes are coupled until $i_1\geq i_0$. If after time
$i_1$, $Y_{i_1+1}^{(a_2)}(\omega)$ becomes positive, then we again
arrive at the inequality $Y_{i_1+1}^{(a_1)}(\omega)\leq
Y_{i_1+1}^{(a_2)}(\omega)$, and so on.

Taking into account that according to the definition $X_{n}^{(0)}$
coincides with $X_n$, we obtain the inequality
\begin{equation}\label{eq6-4}
Y_{n}^{(a)}(\omega)\leq X_{n}(\omega)\leq X_{n}^{(a)}(\omega),
\end{equation}
being correct for all $\omega\in\Omega$ and all $n\geq 0$ as well
as for any initial population $K$ and any $a$. This inequality is
also tight as $a\to 0$, because according to the definition of the
above sequences, $Y_{n}^{(0)}(\omega)=X_{n}^{(0)}(\omega)$ for all
$n$.

Let us now prove the above properties $(i)-(iv)$. Find the limit
in probability of $\frac{Y_n^{(a)}}{K}$ as $K\to\infty$. Notice
first, that according to \eqref{eq2} $\frac{X_1^{(a)}}{K}$
converges to $\max\{a,\frak{m}\}$ in probability, and according to
Wald's equation \cite{Feller 1966}, p.384,
$\frac{\mathrm{E}X_1^{(a)}}{K}$ converge to the same limit
$\max\{a,\frak{m}\}$. Therefore, $\frac{Y_1^{(a)}}{K}$ converges
to $\max\{a,\frak{m}\}-a=\max\{0,\frak{m}-a\}$ in probability, and
$\frac{\mathrm{E}Y_1^{(a)}}{K}$ converges to
$\max\{0,\frak{m}-a\}$. Now, assuming that for some $k$ it is
already proved that $\frac{Y_k^{(a)}}{K}$ converges to
$\max\{0,\frak{m}^k-a\}$ in probability and
$\frac{\mathrm{E}Y_k^{(a)}}{K}$ converges to
$\max\{0,\frak{m}^k-a\}$, by induction we have as follows. If
$\frak{m}^k\leq a$ then $\frac{Y_k^{(a)}}{K}$ converges to 0 in
probability and $\frac{\mathrm{E}Y_k^{(a)}}{K}$ converges to 0,
and consequently,
\begin{equation*}
\begin{aligned}
\mathrm{E}I_k^{(a)} &=\mathrm{P}
\left\{\frac{1}{K}\sum_{j=1}^{Y_k^{(a)}}\xi_{k,j}^\prime>\frac{1}{K}\sum_{j=1}^{\lfloor
aK\rfloor}\Big( 1-\xi_{k,j}\Big)\right\}\to 0.
\end{aligned}
\end{equation*}
The last is true because
\begin{equation*}\frac{1}{K}\mathrm{E}\sum_{j=1}^{Y_k^{(a)}}\xi_{k,j}^\prime=
\frac{1}{K}\mathrm{E}\sum_{j=1}^{Y_k^{(a)}}\xi_{k,j}=\frac{\frak{m}\mathrm{E}
Y_k^{(a)}}{K}\to 0.
\end{equation*}

Therefore, according to \eqref{eq4}
$\frac{\mathrm{E}Y_{k+1}^{(a)}}{K}$ vanishes, and
$\frac{Y_{k+1}^{(a)}}{K}$ vanishes in probability. Therefore, the
assumption $\frak{m}^k\leq a$ is not the case. Hence, assuming
that $\frac{Y_k^{(a)}}{K}$ converges to $\frak{m}^k-a$ in
probability, where $\frak{m}^k>a$, we have the following:
\begin{equation*}
\begin{aligned}
\lim_{K\to\infty}\frac{\mathrm{E}Y_{k+1}^{(a)}}{K} &=\max\{0,
(\frak{m}^k-a)\frak{m}-a(1-\frak{m})\}\\
&=\max\{0,\frak{m}^{k+1}-a\}.
\end{aligned}
\end{equation*}
Thus, as $K\to\infty$, $\frac{Y_n^{(a)}}{K}$ converges to
$\max\{0, \frak{m}^n-a\}$ in probability, and $(i)$ is proved.
Notice, that $(ii)$, $(iii)$ and $(iv)$ follow together with
$(i)$. All these claims are closely related, and their proof is
similar.

Notice also, that the convergence of $\frac{Y_n^{(a)}}{K}$ to
$\max\{0,\frak{m}^n-a\}$ in probability means that in a suitable
probability space, the sequence $\frac{Y_n^{(a)}(\omega)}{K}$
converges almost surely to $\max\{0,\frak{m}^n-a\}$.

\section{\bf Asymptotic expansions and the proof of Lemma \ref{lem1}}\label{sec-4}\noindent
 Pathwise
inequalities  \eqref{eq6-4} and $\frac{Y_n^{(a)}}{K}\leq
\frac{X_n}{K}$ hold for any initial size $K$ and any $a$.
Therefore the appropriate normalized sequences
$\frac{Y_n^{(a)}}{K}$ and $\frac{X_n^{(a)}}{K}$ converge to the
same limit in probability as $K\to\infty$. If there exists the
limit in distribution of $\frac{X_n-\mathrm{E}X_n}{\sqrt{K}}$ as
$K\to\infty$, then because of the equality
$\frac{Y_n^{(a)}-\mathrm{E}Y_n^{(a)}}{\sqrt{K}} =
\frac{X_n^{(a)}-\mathrm{E}X_n^{(a)}}{\sqrt{K}}$, and the
inequality $Y_n^{(a)}(\omega)\leq X_n(\omega)$ for all $a\geq0$
and all $\omega\in\Omega$ (see ref. \eqref{eq6-4}), there are also
the limits in distribution of
$\frac{Y_n^{(a)}-\mathrm{E}Y_n^{(a)}}{\sqrt{K}}$ and
$\frac{X_n^{(a)}-\mathrm{E}X_n^{(a)}}{\sqrt{K}}$ as $K\to\infty$
and $a\to 0$ independently. That is, one can let $K\to\infty$
before $a\to 0$, or converse. Notice, that the limiting
distribution of $\frac{X_n-\mathrm{E}X_n}{\sqrt{K}}$ has been
obtained in \cite{Kl-Ner}, and it also follows from asymptotic
expansions obtained in this section.

It follows from the results of Section \ref{sec-3} that, as
$K\to\infty$, $\tau_{a,K}$ converges in probability to
\begin{equation*}
\ell(a)=\min\{l: \frak{m}^l\leq a\}.
\end{equation*}
and hence, in the case where $K$ increases to infinity first,
$\ell(a)$ = $\mathrm{P}^\_\lim_{K\to\infty}$ $\tau_{a,K}$
($\mathrm{P}^\_\lim$ denotes a limit in probability). It is known
(see e.g. Pakes \cite{Pakes (1989)}), that $\frac{\tau_{K}}{\log
K}$ converges to the constant $c=-\frac{1}{\log\frak{m}}$ in
probability. This result of Pakes \cite{Pakes (1989)} can be
proved by different ways. The advantage of the proof given below
is that it remains true for more general models than the usual
Galton-Watson branching process, resulting in the justice of the
results of the paper for general models as well. For instance, one
can reckon that a bisexual Galton-Watson branching process
starting with $K$ mating units is considered, where $\frak{m}$ now
has the meaning of the average reproduction mean per mating unit
(see Bruss \cite{Bruss 1984}). For the relevant result related to
the $\phi$-branching processes see Bruss \cite{Bruss 1980},
Theorem 1.

For large $X_0=K$ we have as follows:
\begin{equation}\label{4.1}
\begin{aligned}
\tau_K:&=\inf\{t\in\mathbb{N}: X_t=0\}\\ &=\inf\{t\in\mathbb{N}:
X_t<1\}\\
&=\inf\left\{t\in\mathbb{N}: \frac{X_t}{X_0}<\frac{1}{K}\right\}\\
&=\inf\left\{t\in\mathbb{N}:
\prod_{n=1}^t\frac{X_n}{X_{n-1}}<\frac{1}{K}\right\}.
\end{aligned}
\end{equation}
Now note that, as $K\to\infty$, each fraction
$\frac{X_n}{X_{n-1}}$ converge to $\frak{m}$ in probability.
Indeed,
\begin{equation}\label{4.1*}
\frac{X_n}{X_{n-1}}=\frac{X_n}{K}\cdot\frac{K}{X_{n-1}}.
\end{equation}
According to \cite{Kl-Ner}, $\frac{X_n}{K}\to\frak{m}^n$ in
probability as $K\to\infty$. Therefore, the fraction \eqref{4.1*}
converges to $\frac{\frak{m}^n}{\frak{m}^{n-1}}=\frak{m}$ in
probability for any $n$.

On the other hand, by virtue of Wald's identity \cite{Feller
1966}, p.384 we obtain:
\begin{equation}\label{4.1**}
\begin{aligned}
\frac{\mathrm{E}X_n}{\mathrm{E}X_{n-1}}&=\frac{\mathrm{E}\sum_{j=1}^{X_{n-1}}\xi_{n-1,j}}{\mathrm{E}X_{n-1}}
=\frac{\frak{m}\mathrm{E}X_{n-1}}{\mathrm{E}X_{n-1}}=\frak{m}.
\end{aligned}
\end{equation}
So, according to \eqref{4.1*} and \eqref{4.1**}, the limit in
probability of the fraction $\frac{X_n}{X_{n-1}}$ as $K\to\infty$
and the fraction of the corresponding expectations
$\frac{\mathrm{E}X_n}{\mathrm{E}X_{n-1}}$ are the same.

From \eqref{4.1**} we therefore obtain:
\begin{equation}\label{4.2}
\begin{aligned}
\mathrm{E}\left(\frac{X_{t}}{X_0}\right)=\frac{1}{K}\mathrm{E}X_t&=\prod_{n=1}^{t}\frac{\mathrm{E}X_n}{\mathrm{E}X_{n-1}}\\
&=\frak{m}^{t}.
\end{aligned}
\end{equation}
So,
$$\lim_{K\to\infty}K\mathrm{E}\frak{m}^{\tau_K}=1.$$
Similarly to \eqref{4.2}, we also have that
$\left(\frac{X_{t}}{X_0}\right)$ converges to $\frak{m}^{t}$ in
probability as $K\to\infty$ for any integer $t$. So, from
\eqref{4.1} and \eqref{4.2} we have the similar limit as in the
case of expectations, i.e. for any positive $\epsilon$ and
$\delta$ there exists integer $K$ large enough  such that
$\mathrm{P}\{|K\frak{m}^{\tau_K}-1|>\delta\}<\epsilon$, i.e.
 $K\frak{m}^{\tau_K}\to1$
in probability as $K\to\infty$, and thus $\frac{\tau_K}{\log K}$
converges to $-\frac{1}{\log\frak{m}}$ in probability.

From \eqref{eq4} for $Y_{j+1}^{(a)}$, $j=0,1,\ldots$, we obtain
the following equations:
\begin{equation}\label{eq7}
\begin{aligned}
 Y_{j+1}^{(a)}-\frak{m}I_j^{(a)}Y_{j}^{(a)}
  &=\frak{S}\sqrt{I_j^{(a)}K}\sum_{i=1}^{Y_j^{(a)}}\frac{\xi_{j,i}-\frak{m}}
 {\frak{S}\sqrt{I_j^{(a)}K}}-I_j^{(a)}\sum_{i=1}^{\lfloor aK\rfloor}(1-\xi_{j,i})
\end{aligned}
\end{equation}
where ${0}\cdot{\infty}$ is assumed to be 0.

Assuming that $K$ increases to infinity, and dividing both sides
of \eqref{eq7} by large parameter $\frak{S}\sqrt{K}$ we have the
following expansions
\begin{equation}\label{eq8}
\begin{aligned}
&\frac{Y_{j+1}^{(a)}-\frak{m}\chi_{j+1}Y_{j}^{(a)}}{\frak{S}\sqrt{K}}\\
  &\asymp\chi_{j+1}\sum_{i=1}^{Y_j^{(a)}}\frac{\xi_{j,i}-\frak{m}}
 {\frak{S}\sqrt{\chi_{j+1}K}}- \chi_{j+1}\frac{\sqrt{K}}{\frak{S}K}
\sum_{i=1}^{\lfloor aK\rfloor}(1-\xi_{j,i})
\end{aligned}
\end{equation}
or
\begin{equation}\label{eq8-1}
\begin{aligned}
&\chi_{j+1}\cdot
\frac{Y_{j+1}^{(a)}-\frak{m}Y_{j}^{(a)}}{\frak{S}\sqrt{K}}\\
  &\asymp\chi_{j+1}\left(
  \sum_{i=1}^{Y_j^{(a)}}\frac{\xi_{j,i}-\frak{m}}
 {\frak{S}\sqrt{K}}-\frac{a\sqrt{K}}{\frak{S}}\Big(1-\frak{m}\Big)\right).
\end{aligned}
\end{equation}
For $j=0,1,\ldots,\ell(a)-2$, $\ell(a)\geq2$, one can remove the
term $\chi_{j+1}$ from the both sides of \eqref{eq8-1}.

Therefore, for $j=0,1,\ldots,\ell(a)-2$, the left-hand side of
\eqref{eq8-1} can be transformed as follows:
\begin{equation}\label{eq9}
\begin{aligned}
 &\frac{Y_{j+1}^{(a)}-\frak{m}Y_{j}^{(a)}}{\frak{S}\sqrt{K}}\\
  &=\frac{Y_{j+1}^{(a)}-\frak{m}\mathrm{E}Y_{j}^{(a)}
+\frak{m}\mathrm{E}Y_{j}^{(a)}-\frak{m}Y_{j}^{(a)}}{\frak{S}\sqrt{K}}\\
&=\frac{Y_{j+1}^{(a)}-\mathrm{E}Y_{j+1}^{(a)}}{\frak{S}\sqrt{K}}
-\frak{m}\frac{Y_{j}^{(a)}-\mathrm{E}Y_{j}^{(a)}}{\frak{S}\sqrt{K}}
-\frac{a\sqrt{K}}{\frak{S}}\Big(1-\frak{m}\Big)\\
&\asymp\theta_{j+1}^{(a)}-\frak{m}\theta_j^{(a)}-\frac{a\sqrt{K}}
{\frak{S}}\Big(1-\frak{m}\Big),
\end{aligned}
\end{equation}
 where $\{\theta_{j}^{(a)}\}$ is a Gaussian sequence.  (The  values
 of the parameters  of this Gaussian sequence are not
discussed here.)

In turn, for $j=0,1,\ldots,\ell(a)-2$, $\ell(a)\geq 2$, the
right-hand side of \eqref{eq8-1} is transformed as
\begin{equation}\label{eq9-1}
\begin{aligned}
\sum_{i=1}^{Y_j^{(a)}}\frac{\xi_{j,i}-\frak{m}}
 {\frak{S}\sqrt{K}}&-\frac{\sqrt{K}}{\frak{S}K}
 \sum_{i=1}^{\lfloor aK\rfloor}(1-\xi_{j,i})\\
 &=
 \sqrt{\frac{Y_j^{(a)}}{K}}\frac{\xi_{j,i}-\frak{m}}
 {\frak{S}\sqrt{Y_j^{(a)}}}-\frac{a\sqrt{K}}{\frak{S}}\Big(1-\frak{m}\Big)\\
 &=\zeta_j\sqrt{\frac{Y_j^{(a)}}{K}}-
 \frac{a\sqrt{K}}{\frak{S}}\Big(1-\frak{m}\Big),
\end{aligned}
\end{equation}
where $\{\zeta_j\}$ is a sequence of independent standard normally
distributed random variables.

Therefore for $j=0,1,\ldots,\ell(a)-2$ from \eqref{eq9} and
\eqref{eq9-1} we have:
\begin{equation}\label{eq10-1}
\theta_{j+1}^{(a)}-\frak{m}\theta_j^{(a)}
\asymp\zeta_j\sqrt{y_j^{(a)}},
\end{equation}
where $y_j^{(a)}=\frac{Y_j^{(a)}}{K}$.

The analysis of \eqref{eq10-1} is standard. According to the
definition $\mathrm{E}\theta_j^{(a)}=0$. Therefore, rewriting
\eqref{eq10-1} as
\begin{equation}\label{eq11-1}
\theta_{j+1}^{(a)}\asymp\frak{m}\theta_j^{(a)}
+\zeta_j\sqrt{y_j^{(a)}},
\end{equation}
we obtain:
\begin{equation}\label{eq12-1}
\mathrm{E}\left(\theta_{j+1}^{(a)}\right)^2=\lim_{K\to\infty}\mathrm{E}
\left(\frak{m}\theta_j^{(a)}
+\zeta_j\sqrt{\frac{Y_j^{(a)}}{K}}\right)^2.
\end{equation}
Taking into account that
$\mathrm{E}\left(\zeta_j\sqrt{y_j^{(a)}}\right)^2$=$\mathrm{E}y_j^{(a)}\to
 \frak{m}^j-a$, as $K\to\infty$, we obtain
\begin{equation}\label{eq12-2}
\var(\theta_{j+1}^{(a)})=\frak{m}^2\var(\theta_j^{(a)})+
\frak{m}^j-a.
\end{equation}
Next,
 from \eqref{eq10-1} we have:
\begin{equation*}
\begin{aligned}
\cov(\theta_{j}^{(a)}, \theta_{j+1}^{(a)})
&=\mathrm{E}\theta_{j}^{(a)}\theta_{j+1}^{(a)}\\ &=
\lim_{K\to\infty}\mathrm{E}\left(\frak{m}\theta_j^{(a)}\theta_j^{(a)}
+\zeta_j\theta_j^{(a)}\sqrt{\frac{Y_j^{(a)}}{K}}\right)\\
&= \frak{m}\var(\theta_{j}^{(a)})+ \frak{m}^j-a,
\end{aligned}
\end{equation*}
and it is easy to show by induction
\begin{equation*}
\begin{aligned}
&\cov(\theta_{j}^{(a)},
\theta_{j+n}^{(a)})=\frak{m}^n\var(\theta_j^{(a)})+\sum_{i=1}^n\frak{m}^{i-1}
\mathrm{E}y_{j+n-i}^{(a)},\\
&\mathrm{E}y_{j+n-i}^{(a)}=\frak{m}^{j+n-i}-a,
\end{aligned}
\end{equation*}
where $j+n\leq\ell(a)-1$. Assuming now that $a\to 0$, we obtain
the convergence of the sequence
\begin{equation*}
\left\{\frac{X_1-K\frak{m}}{\frak{S}\sqrt{K}},
\frac{X_2-K\frak{m}^2}{\frak{S}\sqrt{K}},\ldots\right\}
\end{equation*}
to the Gaussian process $\{\theta_1, \theta_2,\ldots\}$ with mean
0 and covariance function
\begin{equation}\label{eq-cov}
\begin{aligned}
&\cov(\theta_{j},
\theta_{j+n})=\frak{m}^n\var(\theta_j)+n\frak{m}^{j+n-1},\\
&\var(\theta_{j+1})=\frak{m}^2\var(\theta_j)+\frak{m}^j, \
\var(\theta_1)=1.
\end{aligned}
\end{equation}
This implies the statement of Lemma \ref{lem1}.

\section{\bf Proof of Theorem \ref{th-1}}\label{sec-proof}
\noindent Let us now study equation \eqref{eq10-1} more carefully.
Let $u_1$ and $u_2$ be two real numbers, $0<u_1<u_2<1$. Assume
that $K$ is so large that  the probability
$\mathrm{P}\{|\tau_{a,K}-\ell(a,K)|>\epsilon\}$ is negligible
($\epsilon>0$ is an arbitrary fixed value, $K$ is large enough),
where $\ell(a,K)$ is a (not random) integer number. Such a number
does always exist for any given $a$ since, as $K\to\infty$,
$\tau_{a,K}$ converges to $\ell(a)$ in probability.

For large $K$ we have the following two expansions:
\begin{equation}\label{eq10-2}
\theta_{\lfloor
u_1\tau_{a,K}\rfloor+1}^{(a)}-\frak{m}\theta_{\lfloor
u_1\tau_{a,K}\rfloor}^{(a)} \asymp\zeta_{\lfloor
u_1\tau_{a,K}\rfloor}\sqrt{\frac{Y_{\lfloor
u_1\tau_{a,K}\rfloor}^{(a)}}{K}},
\end{equation}
\begin{equation}\label{eq10-3}
\theta_{\lfloor
u_2\tau_{a,K}\rfloor+1}^{(a)}-\frak{m}\theta_{\lfloor
u_2\tau_{a,K}\rfloor}^{(a)} \asymp\zeta_{\lfloor
u_2\tau_{a,K}\rfloor}\sqrt{\frac{Y_{\lfloor
u_2\tau_{a,K}\rfloor}^{(a)}}{K}},
\end{equation}
where $y_{\lfloor u_i\tau_{a,K}\rfloor}^{(a)}$ in the right-hand
side of equations \eqref{eq10-2} and \eqref{eq10-3}, $i=1,2$, are
correspondingly replaced by $\frac{Y_{\lfloor
u_i\tau_{a,K}\rfloor}^{(a)}}{K}$. It is worth noting as follows.
Relations \eqref{eq10-2} and \eqref{eq10-3} are written in the
form of an asymptotic expansion. The left-hand sides of these
expansions are Gaussian martingale-differences, while the
right-hand sides are the expressions with large parameter $K$.
Since the probability
$\mathrm{P}\{|\tau_{a,K}-\ell(a,K)|>\epsilon\}$ is negligible
($\epsilon>0$ is an arbitrary fixed value, $K$ is large enough),
the expansion with the given right-hand side is correct. From
\eqref{eq10-2} and \eqref{eq10-3} we obtain as follows:
\begin{equation*}
\begin{aligned}
&Y_{\lfloor u_1\tau_{a,K}\rfloor}^{(a)} \left(\theta_{\lfloor
u_2\tau_{a,K}\rfloor+1}-\frak{m}\theta_{\lfloor
u_2\tau_{a,K}\rfloor}\right)^{2}\zeta_{\lfloor
u_1\tau_{a,K}\rfloor}^{2}\\ &\asymp Y_{\lfloor
u_2\tau_{a,K}\rfloor}^{(a)} \left(\theta_{\lfloor
u_1\tau_{a,K}\rfloor+1}-\frak{m}\theta_{\lfloor
u_1\tau_{a,K}\rfloor}\right)^{2}\zeta_{\lfloor
u_2\tau_{a,K}\rfloor}^{2}
\end{aligned}
\end{equation*}
and for any continuous function $f(\bullet)$
\begin{equation}\label{eq10-6-Y-1}
\begin{aligned}
&f\left[Y_{\lfloor u_1\tau_{a,K}\rfloor}^{(a)}
\left(\theta_{\lfloor
u_2\tau_{a,K}\rfloor+1}-\frak{m}\theta_{\lfloor
u_2\tau_{a,K}\rfloor}\right)^{2}\zeta_{\lfloor
u_1\tau_{a,K}\rfloor}^{2}\right]\\ &\asymp f\left[Y_{\lfloor
u_2\tau_{a,K}\rfloor}^{(a)} \left(\theta_{\lfloor
u_1\tau_{a,K}\rfloor+1}-\frak{m}\theta_{\lfloor
u_1\tau_{a,K}\rfloor}\right)^{2}\zeta_{\lfloor
u_2\tau_{a,K}\rfloor}^{2}\right]
\end{aligned}
\end{equation}
For example, from \eqref{eq10-6-Y-1} we obtain:
\begin{equation}\label{eq10-6-Y-1e}
\begin{aligned}
&\left[Y_{\lfloor u_1\tau_{a,K}\rfloor}^{(a)}
\left(\theta_{\lfloor
u_2\tau_{a,K}\rfloor+1}-\frak{m}\theta_{\lfloor
u_2\tau_{a,K}\rfloor}\right)^{2}\zeta_{\lfloor
u_1\tau_{a,K}\rfloor}^{2}\right]^l\\ &\asymp \left[Y_{\lfloor
u_2\tau_{a,K}\rfloor}^{(a)} \left(\theta_{\lfloor
u_1\tau_{a,K}\rfloor+1}-\frak{m}\theta_{\lfloor
u_1\tau_{a,K}\rfloor}\right)^{2}\zeta_{\lfloor
u_2\tau_{a,K}\rfloor}^{2}\right]^l
\end{aligned}
\end{equation}
Now estimate the conditional expectation
$\mathrm{E}\left\{\left(Y_{\lfloor
u_1\tau_{a,K}\rfloor}^{(a)}\right)^l~\Big|~Y_{\lfloor
u_2\tau_{a,K}\rfloor}^{(a)}\right\}$. For brevity let us introduce
a random vector
\begin{equation*}
\mathbf{Z}_{u_1,u_2,\tau_{a,K}}=\left\{\theta_{\lfloor
u_1\tau_{a,K}\rfloor}, \ \theta_{\lfloor u_2\tau_{a,K}\rfloor}, \
\zeta_{\lfloor u_1\tau_{a,K}\rfloor}, \ \zeta_{\lfloor
u_2\tau_{a,K}\rfloor}\right\}.
\end{equation*}
We have
\begin{equation}\label{eq10-6-Y-1a}
\begin{aligned}
&\mathrm{E}\left\{Y_{\lfloor
u_i\tau_{a,K}\rfloor}^{(a)}~\Big|~\mathbf{Z}_{u_1,u_2,\tau_{a,K}}\right\}\\
&=\mathrm{E}\left\{\mathrm{E}\left(Y_{\lfloor
u_i\tau_{a,K}\rfloor}^{(a)}~\Big|~\mathbf{Z}_{u_1,u_2,\tau_{a,K}},
\tau_{a,K}\right)~\Big|~\tau_{a,K}\right\}\\
&=\mathrm{E}\left\{\mathrm{E}\left(Y_{\lfloor
u_i\tau_{a,K}\rfloor}^{(a)}~\Big|~\mathbf{Z}_{u_1,u_2,\tau_{a,K}}\right)~\Big|~\tau_{a,K}
\right\}\\
&=\mathrm{E}Y_{\lfloor u_i\tau_{a,K}\rfloor}^{(a)}, \ \ i=1,2.
\end{aligned}
\end{equation}
The last equality of the right-hand side of \eqref{eq10-6-Y-1a} is
a consequence of conditional independence of $Y_{\lfloor
u_i\tau_{a,K}\rfloor}^{(a)}$ and
$\mathbf{Z}_{u_1,u_2,\tau_{a,K}}$, that is for any given event
$\{\tau_{a,K}=k\}$, the random variable $Y_{\lfloor
u_ik\rfloor}^{(a)}$ and random vector $\mathbf{Z}_{u_1,u_2,k}$ are
independent. \eqref{eq10-6-Y-1a} holds true also in the case of
$a=0$ that will be discussed later.

Next, using the notation $\tau_K=\tau_{0,K}$ let us prove that
\begin{equation}\label{eq11-1-2i}
\cov\left(\theta_{\lfloor
u_1\tau_{K}\rfloor+1}-\frak{m}\theta_{\lfloor u_1\tau_{K}\rfloor},
\theta_{\lfloor u_2\tau_{K}\rfloor+1}-\frak{m}\theta_{\lfloor
u_2\tau_{K}\rfloor}\right)\to 0
\end{equation}
as $K\to\infty$.

Notice, first (see relation \eqref{eq-cov}) that
$\cov(\theta_{j},\theta_{j+n})$ vanishes as $n\to\infty$.
Consequently, by the total expectation formula,
\begin{equation}\label{eq11-1-2j}
\cov(\theta_{\lfloor u_1\tau_K\rfloor},\theta_{\lfloor
u_1\tau_K\rfloor+n})=\mathrm{E}\left(\cov(\theta_{\lfloor
u_1\tau_K\rfloor},\theta_{\lfloor
u_1\tau_K\rfloor+n}~|~\tau_K)\right)
\end{equation}
vanishes as $n\to\infty$, where here in relation \eqref{eq11-1-2j}
and later the notation for $\cov(\theta_{\lfloor
u_1\tau_K\rfloor}, \theta_{\lfloor u_1\tau_K\rfloor+n}~|~\tau_K)$
or another similar notation means the conditional covariance.
Taking into account that, as $K\to\infty$, $\tau_K$ increases to
infinity in probability and $u_2-u_1>0$, the difference $\lfloor
u_2\tau_K\rfloor-\lfloor u_1\tau_K\rfloor$ increases to infinity
in probability too. Hence, by virtue of \eqref{eq11-1-2j} one can
conclude that $\cov(\theta_{\lfloor
u_1\tau_K\rfloor},\theta_{\lfloor u_2\tau_K\rfloor})$ vanishes as
$K\to\infty$. Therefore, as $K\to\infty$,
$\mathrm{E}\theta_{\lfloor u_1\tau_K\rfloor}\theta_{\lfloor
u_2\tau_K\rfloor}$ is asymptotically equal to
$\mathrm{E}\theta_{\lfloor
u_1\tau_K\rfloor}\mathrm{E}\theta_{\lfloor u_2\tau_K\rfloor}$, and
\eqref{eq11-1-2i} follows. In addition to \eqref{eq10-6-Y-1a} and
\eqref{eq11-1-2i} we have also the following. Since the sequence
$\{\zeta_j\}$ consists of independent standard normally
distributed random variables, then as $K\to\infty$
\begin{equation}\label{eq11-1-2k}
\cov\left(\zeta_{\lfloor u_1\tau_K\rfloor},\zeta_{\lfloor
u_2\tau_K\rfloor}\right)\to 0.
\end{equation}
This is because $\cov\left(\zeta_{\lfloor
u_1\tau_K\rfloor},\zeta_{\lfloor
u_2\tau_K\rfloor}~|~\tau_K\right)=\mathrm{I}\{\lfloor
u_1\tau_K\rfloor=\lfloor u_2\tau_K\rfloor\}$, and the last
vanishes in probability as $K\to\infty$.

Assuming  that $a$ vanishes we need a stronger assumption than
above. Specifically, we assume that $K$ is so large that  the
probability
\begin{equation*}
\mathrm{P}\left\{\Big|\frac{{\tau_{a,K}-\ell(a,K)}}{\frac{\log
K}{\log\frak{m}}}\Big|>\epsilon\right\}
\end{equation*}
 is negligible for all $0\leq a<a_0$ ($\epsilon>0$ is an arbitrary
 fixed value, $K$ is large enough), where $a_0<1$ is some fixed small
number. Such a large number $K$ does always exist, since as
$K\to\infty$ and $a$ vanishing, $\frac{\tau_{a,K}}{\log K}$
converges to $-\frac{1}{\log\frak{m}}$ in probability. Then,
letting $a\to 0$ in \eqref{eq10-6-Y-1e} in view of pathwise
inequalities \eqref{eq6-4} and $\frac{Y_n^{(a)}}{K}\leq
\frac{X_n}{K}$ we have
\begin{equation}\label{eq11-1-2l}
\begin{aligned}
&X_{\lfloor u_1\tau_{K}\rfloor}^l \left(\theta_{\lfloor
u_2\tau_{K}\rfloor+1}-\frak{m}\theta_{\lfloor
u_2\tau_{K}\rfloor}\right)^{2l}\zeta_{\lfloor
u_1\tau_{K}\rfloor}^{2l}\\ &\asymp X_{\lfloor
u_2\tau_{K}\rfloor}^l \left(\theta_{\lfloor
u_1\tau_{K}\rfloor+1}-\frak{m}\theta_{\lfloor
u_1\tau_{K}\rfloor}\right)^{2l}\zeta_{\lfloor
u_2\tau_{K}\rfloor}^{2l}.
\end{aligned}
\end{equation}
Taking into account \eqref{eq10-6-Y-1a}, \eqref{eq11-1-2i} and
\eqref{eq11-1-2k} and conditional independency of $X_{\lfloor
u_1\tau_{K}\rfloor}$, $\left(\theta_{\lfloor
u_1\tau_{K}\rfloor+1}-\frak{m}\theta_{\lfloor
u_1\tau_{K}\rfloor}\right)$ and $\zeta_{\lfloor
u_2\tau_{K}\rfloor}$, and passing to the appropriate conditional
expectations, from \eqref{eq11-1-2l} we obtain:
\begin{equation}\label{eq10-6-X-3}
\begin{aligned}
\mathrm{E}\left\{X_{\lfloor u_1\tau_{K}\rfloor}^l~\Big|~X_{\lfloor
u_2\tau_{K}\rfloor}\right\} &=
\mathrm{E}\left[\mathrm{E}\left\{X_{\lfloor
u_1\tau_{K}\rfloor}^l~\Big|~X_{\lfloor u_2\tau_{K}\rfloor},
\tau_K\right\}
~\Big|~\tau_K\right]\\
&\asymp X_{\lfloor
u_2\tau_{K}\rfloor}^l\frac{\mathrm{E}\left(\theta_{\lfloor
u_1\tau_{K}\rfloor+1}-\frak{m}\theta_{\lfloor
u_1\tau_{K}\rfloor}\right)^{2l}} {\mathrm{E}\left(\theta_{\lfloor
u_2\tau_{K}\rfloor+1}-\frak{m}\theta_{\lfloor
u_2\tau_{K}\rfloor}\right)^{2l}}.
\end{aligned}
\end{equation}

Thus, to this end our task is to determine the asymptotic of
\begin{equation*}
\mathrm{E}\left(\theta_{\lfloor
u_i\tau_{K}\rfloor+1}-\frak{m}\theta_{\lfloor
u_i\tau_{K}\rfloor}\right)^{2l}, \ \ i=1,2,
\end{equation*}
for large $K$. Returning to basic equations \eqref{eq10-2} and
\eqref{eq10-3}, we have
\begin{equation}\label{eq11-1-3}
\begin{aligned}
&\mathrm{E}\left(\theta_{\lfloor u_i\tau_{a,K}\rfloor+1}^{(a)}
-\frak{m}\theta_{\lfloor u_i\tau_{a,K}\rfloor}^{(a)}\right)^{2l}
\asymp\mathrm{E}\left(\zeta_{\lfloor
u_i\tau_{a,K}\rfloor}\sqrt{\frac{Y_{\lfloor
u_i\tau_{a,K}\rfloor}^{(a)}}{K}}\right)^{2l}\\
&=\mathrm{E}\left\{\mathrm{E}\left[\left(\zeta_{\lfloor
u_i\tau_{a,K}\rfloor}\sqrt{\frac{Y_{\lfloor
u_i\tau_{a,K}\rfloor}^{(a)}}{K}}\right)^{2l}~\Big|~\tau_{a,K}\right]\right\}\\
&=\mathrm{E}\left(y_{\lfloor u_i\tau_{a,K}\rfloor}^{(a)}\right)^l
\mathrm{E}(\zeta_{\lfloor u_i\tau_{a,K}\rfloor})^{2l},
\end{aligned}
\end{equation}
where $\zeta_{\lfloor u_i\tau_{a,K}\rfloor}$, $i=1,2$, are
standard normally distributed random variables. As $a$ vanishes,
from \eqref{eq11-1-3} we obtain
\begin{equation}\label{eq11-1-4}
\begin{aligned}
\mathrm{E}\left(\theta_{\lfloor
u_i\tau_{K}\rfloor+1}-\frak{m}\theta_{\lfloor
u_i\tau_{K}\rfloor}\right)^{2l}\asymp\mathrm{E}\frak{m}^{l\lfloor
u_i\tau_{K}\rfloor}\mathrm{E}(\zeta_{\lfloor
u_i\tau_{K}\rfloor})^{2l}.
\end{aligned}
\end{equation}
Therefore, \eqref{eq10-6-X-3} can be rewritten
\begin{equation*}
\begin{aligned}
\mathrm{E}\left\{X_{\lfloor u_1\tau_{K}\rfloor}^l~\Big|~X_{\lfloor
u_2\tau_{K}\rfloor}\right\}\asymp X_{\lfloor
u_2\tau_{K}\rfloor}^l\mathrm{E}\frak{m}^{l(\lfloor
u_1\tau_{K}\rfloor-\lfloor u_2\tau_{K}\rfloor)}.
\end{aligned}
\end{equation*}
This proves the first equation of \eqref{eq-th-1}. The proof of
the second equation of \eqref{eq-th-1} is similar.

Consider basic equation \eqref{eq10-2} again, rewriting it as
follows:
\begin{equation}\label{eq10-2end}
Y_{\lfloor u_1\tau_{a,K}\rfloor}^{(a)}\zeta_{\lfloor
u_1\tau_{a,K}\rfloor}^2\asymp K\left(\theta_{\lfloor
u_1\tau_{a,K}\rfloor+1}^{(a)}-\frak{m}\theta_{\lfloor
u_1\tau_{a,K}\rfloor}^{(a)}\right)^2.
\end{equation}
Assuming that as $a$ vanishes we have:
\begin{equation}\label{eq10-2end-1}
X_{\lfloor u_1\tau_{K}\rfloor}\zeta_{\lfloor
u_1\tau_{K}\rfloor}^2\asymp K\left(\theta_{\lfloor
u_1\tau_{K}\rfloor+1}-\frak{m}\theta_{\lfloor
u_1\tau_{K}\rfloor}\right)^2.
\end{equation}
Therefore, taking into account that $X_{\lfloor
u_1\tau_{K}\rfloor}$ and $\zeta_{\lfloor u_1\tau_{K}\rfloor}$ are
conditionally independent, from \eqref{eq10-2end-1} we obtain:
\begin{equation*}
\begin{aligned}
\mathrm{E}\left\{X_{\lfloor
u_1\tau_{K}\rfloor}^l~\Big|~\tau_{K}\right\} &\asymp
K^l\frac{\mathrm{E}\left\{ \left(\theta_{\lfloor
u_1\tau_{K}\rfloor+1}-\frak{m}\theta_{\lfloor
u_1\tau_{K}\rfloor}\right)^{2l}~|~\tau_K\right\}}{\mathrm{E}(\zeta_{\lfloor
u_1\tau_{a,K}\rfloor}~|~\tau_K)^{2l}}\\
&\asymp K^l\frak{m}^{l\lfloor u_1\tau_K\rfloor}.
\end{aligned}
\end{equation*}
\eqref{eq-th-3} is proved.

\section{\bf Discussion}\label{sec-5}
\noindent The aim of this section is to present the main results
in convenient form for application to analysis of real
populations. In this section we also establish a so-called
invariance property.

Let, when $K$ is large, $\epsilon$ be a relatively small (positive
or negative) parameter having the following meaning. The
population size at time $\lfloor u_2\tau_K\rfloor$ is assumed to
be equal to $\lfloor(1+\epsilon)K\frak{m}^{u_2t_K}\rfloor$,
$t_K=-\frac{\log K}{\log\frak{m}}$.

The meaning of this value is the following. The factor
$K\frak{m}^{\lfloor u_2t_K\rfloor}$ is the expected size of the
population at time $\lfloor u_2\tau_K\rfloor$, and the factor
$1+\epsilon$ represents a parameter of relative deviation from the
expected population at that time moment. Then, from Theorem
\ref{th-1} we obtain, that for large $K$
\begin{equation}\label{eq12-5}
\begin{aligned}
&\log\mathrm{E}\left(X_{\lfloor u_1\tau_K\rfloor}^l~|~X_{\lfloor
u_2\tau_K\rfloor}=\lfloor(1+\epsilon)K\frak{m}^{u_2t_K}\rfloor\right)\\
&=l\log(K+K\epsilon)+lu_1 t_K\log\frak{m} +o(1)\\
&=l\log(K+K\epsilon)-lu_1\log K +o(1).
\end{aligned}
\end{equation}
In real computations the term $\log (K+K\epsilon)$ can be replaced
by $\epsilon+\log K$ if $\epsilon$ is sufficiently small.

The result similar to \eqref{eq12-5} can be obtained for the
conditional expectation of \eqref{eq-th-3}. Specifically, for
large $K$ write
\begin{equation}\label{eq12-6}
\tau_K=-\left\lfloor(1+\epsilon)\frac{\log K}{\log
\frak{m}}\right\rfloor.
\end{equation}
\eqref{eq12-6} has the following meaning. As $K\to\infty$, the
fraction $\frac{\tau_K}{\log K}$ converges in probability to
$-\frac{1}{\log\frak{m}}$, and therefore, as $K$ is large, the
factor 1+$\epsilon$ is a parameter for relative deviation from the
expected value of extinction time. Then,
\begin{equation}\label{eq12-7}
\begin{aligned}
&\log\mathrm{E}\left\{X_{\lfloor
u_1\tau_K\rfloor}^l~\Big|~\tau_K=-\left\lfloor(1+\epsilon)\frac{\log
K}{\log\frak{m}}\right\rfloor\right\}\\&=l\log
(K+K\epsilon)-lu_1\log K+o(1).
\end{aligned}
\end{equation}
As we can see the right-hand sides of \eqref{eq12-5} and
\eqref{eq12-7} coincide. That is for any given relative deviation
1+$\epsilon$ the asymptotic conditional expectations are
invariant.

\section*{\bf Acknowledgement}
The author thanks Prof. Peter Jagers for useful conversation and
advice.
The research
of the author is supported by the Australian Research Council,
Grant \# DP0771338.

\bibliographystyle{amsplain}

\end{document}